\documentclass{amsart}
\usepackage{amssymb,latexsym,graphics}
\usepackage[mathscr]{eucal}

\newtheorem{thm}{Theorem}[section]

\newtheorem{lem}[thm]{Lemma}

\theoremstyle{remark}

\theoremstyle{definition}

\newcommand{\Lp}[2]{\left\Vert \, #1 \, \right\Vert_{#2}}
\newcommand{\lp}[2]{\Vert \, #1 \, \Vert_{#2}}

\newcommand{\ret}{\vspace{.4cm}}

\begin{document}

\title[Global Regularity for NLW]
{Global Regularity for General Non-Linear Wave Equations \textbf{I}. 
$(6+1)$ and Higher Dimensions}
\author{}
\address{}
\email{}
\author{Jacob Sterbenz}
\address{Department of Mathematics, Princeton University,
Princeton NJ, 08544}
\email{sterbenz@math.princeton.edu}
\thanks{This work was conducted under NSF grant DMS-0100406.}
\subjclass{}
\keywords{}
\date{}
\dedicatory{}
\commby{}


\begin{abstract}
Following work of Tataru, \cite{Tataru} and \cite{tatwm1}, 
we  solve the division problem
for wave equations with generic quadratic non-linearities in high dimensions.
Specifically, we show that non--linear wave equations 
which can be written as systems involving equations of the form  
$\Box \phi = \phi\,\nabla\phi$ and $\Box \phi = |\nabla\phi|^2$
are well-posed with scattering
in $(6+1)$ and higher dimensions if the Cauchy data are small in the 
scale invariant $\ell^1$ Besov space $\dot{B}^{s_c,1}$.  
This paper is the first in a series of works where we discuss 
the global regularity 
properties of general non-linear wave equations for all dimensions $4\leqslant n$.
\end{abstract}

\maketitle


\section{Introduction}

In this paper, our aim is to give a more or less complete description
of the global regularity properties of generic homogeneous quadratic
non--linear wave equations
on $(6+1)$ and higher dimensional Minkowski space. The equations
we will consider are all of the form:
\begin{equation}\label{generic_system}
	\Box \phi \ = \ \mathcal{N}(\phi,D\phi) \ . 
\end{equation}
Here $\Box := -\partial_t^2 + \Delta_x$ denotes the standard wave operator
on $\mathbb{R}^{n+1}$, and 
$\mathcal{N}$ is a smooth function of $\phi$ and its first partial
derivatives, which we denote by $D\phi$. For all of the nonlinearities we
study here, $\mathcal{N}$ will be assumed 
to be at least quadratic in nature, that is:
\begin{align}
	\mathcal{N}(X,Y) \ &= \ O(|(X,Y)|^2) \ ,  &(X,Y) \ 
	&\sim \ 0 \ . \notag  
\end{align}
The homogeneity condition we require $\mathcal{N}$ to satisfy is that there 
exist a (vector) $\sigma$ such that:
\begin{equation}\label{homog_condition}
	\mathcal{N}(\lambda^\sigma\phi,
	\lambda^{\sigma +1}
	D\phi) \ = \ \lambda^{\sigma + 2}
	\mathcal{N}(\phi,
	D\phi) \ ,
\end{equation}	
where we use multiindex notation for vector $\mathcal{N}$.
The condition \eqref{homog_condition}
implies that solutions to the system \eqref{generic_system} are
invariant (again solutions) if one performs the scale transformations:
\begin{equation}\label{scale_trans}
	\phi(\cdot) \ \rightsquigarrow \ 
	\lambda^\sigma\phi(\lambda\, \cdot) \ .
\end{equation}
The general class of equations which falls under this
description contains virtually all
massless non--linear field theories on Minkowski space, including the
Yang Mills equations (YM), the wave--maps equations (WM), and the 
Maxwell--Dirac equations (MD). We list the schematics for these 
systems respectively as:\\

\begin{align}
	\Box A \ &= \ A\, DA + A^3 \ , \tag{YM} \\ \notag \\
	\Box \phi \ &= \ |D\phi|^2 \ , \tag{WM} \\ \notag \\
	\begin{split}
		\Box u \ &= \ A\, Du \ , \\
		\Box A \ &= \ |Du|^2 \ .
	\end{split}\tag{MD}
\end{align}\\

\noindent The various values of $\sigma$ for these equations are
(respectively) $\sigma=1$, $\sigma=0$, and $\sigma=(\frac{1}{2},1)$.
For a more complete introduction to these equations, see for instance
the work \cite{FKbilin} and \cite{Bmd}. For the purposes of this paper,
will will only be concerned with the structure of these equations at the 
level of the generic schematics (YM)--(MD).\\

The central problem we will be concerned with is that of giving a
precise description of the
regularity assumptions needed in order to guarantee that the
Cauchy problem for the system \eqref{generic_system} is globally well
posed with scattering (GWPS). That is, given initial data:
\begin{align}\label{cauchy_data}
	&\phi(0)\ = \ f  \ , &\partial_t\phi\, (0)\ = \ g \ , 
\end{align}
we wish to describe how much smoothness and decay $(f,g)$ needs to
possess in order for there to exist a unique global solution to the system
\eqref{generic_system} with this given initial data. We also wish to show
that the solutions we construct depend continuously on the initial data, and
are asymptotic to solutions of the linear part of 
\eqref{generic_system}. We will describe shortly in what sense we will
require these notions to hold.\\

Our main motivation here is to be able to prove global well--posedness for non-linear
wave equations of the form \eqref{generic_system}
in a context where the initial data may not be very smooth, and furthermore does
not possess enough decay at space--like infinity to be in $L^2$. 
Also, we would like to understand how this can be done in situations where the
equations being considered contain no special structure in the non-linearity.
For instance, this is of interest in discussing the problem of small data
global well--posedness for the Maxwell--Klein--Gordon and Yang--Mills equations 
with the Lorentz gauge enforced instead of the more regular Coulomb gauge. 
This provides a significant point of departure from 
earlier works on the global existence theory of non--linear wave equations, 
which for the general case requires precise control on the initial data in certain 
weighted Sobolev spaces  (see \cite{Kvect}), or else requires the non--linearity to
have some special algebraic or gauge structure which allows one to exploit 
some null form identities or apply an appropriate renormalization to
the equation being considered (see \cite{tatwm1} and \cite{Taowm}).\\
 
From the point of view of homogeneity, we are lead directly to considerations
of the low regularity properties of the equations \eqref{generic_system} as follows:
By a simple scaling argument\footnote{In conjunction with finite time blowup 
for large data. This phenomena is known to happen for higher dimensional 
equations with derivative non--linearities even in the presence of positive conserved 
quantities.}:
one can see that the most efficient $L^2$ based regularity assumption 
possible on the initial data involves $s_c = \frac{n}{2} - \sigma$ derivatives.
Again, by scale invariance and
looking at unit frequency initial data\footnote{That is initial data sets where
the Fourier transform is supported in the unit frequency annulus $\{\xi \ \big| \
\frac{1}{2} < |\xi| < 2\}$.}, 
one can see that if we are to impose
only an $L^2$ smallness condition on the initial data which contains no 
physical space weights, then
$s_c = \frac{n}{2} - \sigma$ is in fact the \emph{largest} amount of
derivatives we may work with.
This leads us to consider the question of GWPS for initial data
in the homogeneous Besov spaces $\dot{B}^{s_c,p}$, for various values of $p$.
In this work, we will concentrate solely on the case $p=1$. This is the strongest
scale and translation invariant control on the initial data possible, and will be 
crucial for the kind of non--linearities we work with here. In fact, it does not 
seem possible to push any type of global regularity for equations of the
type \eqref{generic_system} which contain derivatives in the non--linearity
down to the scale invariant Sobolev space $\dot H^{s_c} = \dot B^{s_c,2}$ unless
the equations under consideration possess a great deal of special structure in the
non--linearity. This has been done for the wave--maps equations (see \cite{Taowm}) and
more recently for the Maxwell--Klein--Gordon equations in $(6+1)$ and higher dimensions
with the help of the Coulomb gauge (see \cite{RTmkg}). \\

In recent years, there
has been much progress in our understanding of the low regularity local theory
for general non--linear wave equations of the form \eqref{generic_system}.
In the lower dimensional setting, i.e. when $n=2,3,4$, it is known from
counterexamples of Lindblad (see \cite{Lind_counter}) 
that there is ill posedness
for initial data in the Sobolev space $H^{s_0}$, where
$s_0 \leqslant s_c + \frac{5-n}{4}$.
Intimately connected with this phenomena is the failure of certain space--time 
estimates for the linear wave equation known as \emph{Strichartz estimates}.
Specifically, one does not have anything close to an $L^2(L^4)$ estimate in 
these dimensions. Such an estimate obviously plays a crucial role (via 
Duhamel's principle) in the quadratic theory. However, 
using the Strichartz estimates available in these dimensions
along with Picard iteration in certain function spaces, 
one can show that the Lindblad counterexamples are sharp in
there is local well--posedness for
initial data in the spaces $H^s$ when $s_c + \frac{5-n}{4} < s$ 
(see for example \cite{KSel}).\\
 
In the higher dimensional setting, i.e. when the number of spatial dimensions 
is  $n=5$ or greater, one does have access to Strichartz estimates at the level of
$L^2(L^4)$ (see \cite{KT_str}), and it is possible to push the local 
theory down to $H^{s_c + \epsilon}$, where $0 < \epsilon$ is arbitrary
(see \cite{Tataru}).\\

In all dimensions, the single most important factor which determines
the local theory as well as the range of validity for Strichartz estimates
is the existence of free waves which are highly
concentrated along null directions in Minkowski space. These waves,
known as \emph{Knapp counterexamples}, resemble a single beam of light
which remains coherent for a long period of time before dispersing.
For a special class of non--linearities, known as ``null structures'',
interactions between these coherent beams are effectively canceled,
and one gains an improvement in the local theory of equations whose 
nonlinearities have this form (see for example \cite{KMQ0}, \cite{KSel}).\\

In both high and low dimensional settings, the analysis of certain
null structures, specifically non--linearities containing
the $Q_0$ null from\footnote{This is defined by the equation $Q_0(\phi,\psi) = 
\partial_\alpha\phi\, \partial^\alpha\psi$.}, 
has led to the proof that the wave--maps model
equations\footnote{ Not the rough schematic we have listed here, but rather 
equations of the form \newline $\Box\phi = \Gamma(\phi)Q_0(\phi,\phi)$.} 
are well  posed in the scale invariant $\ell^1$ Besov space 
$\dot{B}^{\frac{n}{2},1}$ (see \cite{tatwm1} and \cite{tatwm2}). 
While the proof of this result is quite simple 
for high dimensions, it relies in an essential way on the 
structure of the $Q_0$ null form. In fact, there is no direct way to
extend the proof of this result to include the less regular nonlinearities
of the form $\phi\, \nabla \phi$, or for that matter the
$Q_{ij}$ null forms\footnote{These are defined by $Q_{ij}(\phi,\psi) = 
\partial_i \phi , \partial_j \psi - \partial_j \phi \partial_i \psi$.},
which show up in the equations of gauge field theory.
However, the high dimensional non--linear interaction
of coherent waves is quite weak (e.g. giving the desired range of validity
for Strichartz estimates), and one would expect that it is possible to prove
local well posedness for quadratic equations with initial data
in the scale invariant $\ell^1$
Besov space without resorting to any additional structure in the 
nonlinearity. For $n=5$ dimensions, it may be that this is not quite possible,
although we provide no convincing evidence except for the fact that there is
no obvious way to add over our localized estimates in that dimension in order to
obtain a full set of estimates that works in all of space--time Fourier space. For
$n=6$ and higher dimensions, we will prove that in fact no null structure
is needed for there to be well posedness in $\dot{B}^{s_c,1}$. This leads
to the statement of our main result which is as follows:\\
 
\begin{thm}[Global Well Posedness]\label{WP_theorem}
	Let $6 \leqslant n$ be the number of spatial dimensions. For 
	any of the generic equations listed above: 
	YM, WM, or MD, let $(f,g)$ be a (possibly vector valued) initial
	data set. Let $s_c = \frac{n}{2} - \sigma$ be the corresponding
	$L^2$ scaling exponent.
	Then there exists constants $0 < \epsilon_0,C$ such that if
	\begin{equation}\label{initial_smallness}
		\lp{(f,g)}{\dot{B}^{s_c,1}\times\dot{B}^{s_c-1,1}}
		\ \ \leqslant \ \ \epsilon_0 \ ,
	\end{equation}
	there exits a global solution $\psi$ 
	which satisfies the continuity condition:
	\begin{equation}\label{cont_cond}
		\lp{\psi}{C(\dot{B}^{s_c,1})\cap 
		C^{(1)}(\dot{B}^{s_c-1,1})} \ \ 
		\leqslant \ \ C\lp{(f,g)}{\dot{B}^{s_c,1}\times
		\dot{B}^{s_c-1,1}} \ .
	\end{equation}
	The solution $\psi$ is unique in the following sense: There
	exists a sequence of smooth functions $(f_N,g_N)$ such that: 
	\begin{equation}
		\lim_{N\to \infty}
		\lp{(f,g) - (f_N,g_N)}{\dot{B}^{s_c,1}\times\dot{B}^{s_c-1,1}}
		 \ =  \ 0 \ . \notag
	\end{equation}
	For this sequence of functions, there exists a sequence of unique
	smooth global solutions $\psi_N$ of \eqref{generic_system}
	with this initial data. Furthermore, the $\psi_N$ converge to
	$\psi$ as follows:
	\begin{equation}
		\lim_{N\to \infty}
		\lp{\psi - \psi_N}{ C(\dot{B}^{s_c,1})\cap
		C^{(1)}(\dot{B}^{s_c-1,1}) }  \ = \ 0 \ . \notag
	\end{equation}  
	Also, $\psi$ is the only solution which may be obtained as a 
	limit (in the above sense) of solutions to \eqref{generic_system}
	with regularizations of ${(f,g)}$ as initial data. 
	Finally, $\psi$ retains any extra smoothness inherent in the initial
	data. That is, if $(f,g)$ also has finite $\dot{H}^{s}\times 
	\dot{H}^{s-1}$ norm, for $s_c < s$, then so does 
	$\psi$ at fixed time and one has the following estimate:
	\begin{equation}\label{smooth_cont_cond}
		\lp{\psi}{C(\dot{H}^s)\cap C^{(1)}(\dot{H}^{s-1})} \ \ 
		\leqslant \ \ C\lp{(f,g)}{\dot{H}^s\times \dot{H}^{s-1}} \ .
	\end{equation}
\end{thm}\ret\ret

\noindent In a straightforward way, the function spaces we iterate in allow us
to show the following  scattering result without any extra work:\\ 

\begin{thm}\label{scattering_result}
	Using the same notation as above we have that there exists
	data sets $(f^\pm,g^\pm)$, such that if
	$\psi^\pm$ is the solution to the homogeneous wave equation with
	the corresponding initial data, the following asymptotics hold:
	\begin{align}
		\lim_{t \rightarrow \infty}
		\lp{\psi^+ - \psi}{\dot{B}^{s_c,1}\cap\partial_t 
		\dot{B}^{s_c-1,1}}
		\ &= \ 0 \ , \label{scattering_result1} \\
		 \lim_{t \rightarrow -\infty}
		\lp{\psi^- - \psi}{\dot{B}^{s_c,1}\cap\partial_t 
		\dot{B}^{s_c-1,1}}
		\ &= \ 0 \ . \label{scattering_result2} 
	\end{align}
	Furthermore, the scattering operator retains any additional regularity
	inherent in the initial data. That is, if $(f,g)$ has finite
	$\dot{H}^s$ norm, then so does $(f^\pm,g^\pm)$, and the following
	asymptotics hold:
	\begin{align}
		\lim_{t \rightarrow \infty}
		\lp{\psi^+ - \psi}{\dot{H}^s\cap\partial_t \dot{H}^{s-1}}
		\ &= \ 0 \ , \label{smooth_scattering_result1} \\
		 \lim_{t \rightarrow -\infty}
		\lp{\psi^- - \psi}{\dot{H}^s\cap\partial_t \dot{H}^{s-1}}
		\ &= \ 0 \ . \label{amooth_scattering_result2} 
	\end{align} 
\end{thm}\ret

\ret

\section{Preliminary Notation}

For quantities $A$ and $B$, we denote by $A \lesssim B$ to mean
that $A \leqslant C\cdot B$ for some large constant $C$. The constant
$C$ may change from line to line, but will always remain
fixed for any given instance where this notation appears.
Likewise we use the notation $A\sim B$ to mean that
$\frac{1}{C}\cdot B \leqslant A \leqslant C\cdot B$. We also
use the notation $A \ll B$ to mean that $A \leqslant
\frac{1}{C}\cdot B$ for some
large constant $C$. This is the notation we will use throughout
the paper to break down quantities into the standard cases:
$A\sim B$, or $A \ll B$, or $B \ll A$; and
$A \lesssim B$, or $B \ll A$, without ever discussing
which constants we are using. \\

For a given function of two variables 
$(t,x)\in \bf{R}\times\bf{R}^3$ we write the spatial
and space--time Fourier transform as:
\begin{align}
	\widehat{f}(t,\xi) &= \int e^{-2\pi i \xi\cdot x} \, f(t,x) \ dx
	\ , \notag \\
	\widetilde{f}(\tau,\xi) &= \int e^{-2\pi i (\tau t + \xi\cdot x)} 
	\ f(t,x) \ dt  dx \ . \notag
\end{align}
respectively. At times, we will also write $\mathscr{F}[f] = 
\widetilde{f}$.\\

For a given set of functions of the spatial variable only, we
denote by $W(f,g)$ the solution of the homogeneous wave equation 
with Cauchy data $(f,g)$.
If $F$ is a function on space--time, we will denote by $W(F)$ 
the function $W\left(F(0),\partial_t F\, (0)\right)$.\\
 
Let $E$ denote any fundamental solution to the homogeneous wave equation.
i.e., one has the formula $\Box E = \delta$. We define the standard
Cauchy parametrix for the wave equation by the formula:
\begin{equation}
        \frac{1}{\Box} F \ = \ E*F - W (E*F) \ . \notag
\end{equation}
Explicitly, one has the identity:
\begin{equation}\label{Duhamel_integral}
        \widehat{\frac{1}{\Box}F}\, (t,\xi) \
        = \ - \int_0^t \frac{\sin\left(2\pi
        |\xi|(t-s)\right)}{2\pi |\xi|} \widehat{F}(s,\xi) \ ds \ . 
\end{equation}\\
 
For any function $F$ which is supported away from the light cone in 
Fourier space, we shall use the following notation for division by the
symbol of the wave equation:
\begin{equation}
        \frac{1}{\varXi} F \ = \
        E*F \ . \notag
\end{equation}
Of course, the definition of $\frac{1}{\varXi}$ does not depend on $E$ 
so long as for $F$ is supported away from the light cone; 
for us that will always be the case when we use this notation. 
Explicitly, one has the formula:
\begin{equation}
        \mathcal{F}\left[ \frac{1}{\varXi} F
        \right](\tau,\xi) \
        = \ \frac{1}{4\pi^2(\tau^2 - |\xi|^2)}
	\widetilde{F}(\tau,\xi) \ . \notag
\end{equation}\\

\ret

\section{Multipliers and Function Spaces}

Let $\varphi$ be  a smooth bump function (i.e.
supported on the set $|s| \leqslant 2$ such
that $\varphi = 1$ for $|s| \leqslant 1$). In what follows, it will be
a great convenience for us to assume that $\varphi$ may change its exact
form for two separate instances of the symbol $\varphi$ (even if they
occur on the same line). In this way, we may assume without loss of generality 
that in addition to being smooth, we also have the idempotence identity
$\varphi^2= \varphi$. We shall use this convention for all the cutoff functions
we introduce in the sequel.\\

For $\lambda\in 2^\mathbb{Z}$, 
we denote the dyadic scaling of $\varphi$ by 
$\varphi_\lambda(s) = \varphi(\frac{s}{\lambda})$. The most basic Fourier
localizations we shall use here are with respect to the space-time variable
and the distance from the cone. Accordingly, we form the Littlewood-Paley
type cutoff functions:
\begin{align}
	s_{\lambda}(\tau,\xi) &= \varphi_{2\lambda}(|(\tau,\xi)|)
	- \varphi_{\frac{1}{2}\lambda}(|(\tau,\xi)|) \ , \label{st_cutoff} \\
	c_d(\tau,\xi) &= \varphi_{2d}(|\tau| - |\xi|) - 
	\varphi_{\frac{1}{2} d}(|\tau| - |\xi|) \ . \label{cone_cutoff} 
\end{align}
We now denote the corresponding Fourier multiplier operator via the formulas 
$\ \ \widetilde{S_\lambda u} = s_\lambda \widetilde{u}\ \ $ and
$\ \ \widetilde{C_d u} = c_d \widetilde{u}\ \ $ respectively. We also use a 
multi-subscript notation to denote products of the above operators, e.g.
$\ \ S_{\lambda,d} = S_\lambda C_d\ \ $. We shall use the notation:
\begin{equation}
	S_{\lambda,\bullet \leqslant d} = \sum_{\delta \leqslant d}
	S_{\lambda,\delta} \ , \label{d_cone_nbd}
\end{equation}
to denote cutoff in an $O(d)$ neighborhood of the light cone in Fourier
space. At times it will also be convenient to write 
$\ \ S_{\lambda,d \leqslant \bullet} = S_\lambda - 
S_{\lambda,\bullet < d}\ \ $. We shall also use the notation
$S^\pm_{\lambda,d}$ etc. to denote the multiplier $S_{\lambda,d}$ 
cutoff in the half space $\ \ \pm \tau > 0\ \ $.\\

The other type of Fourier localization which will be central to our analysis
is the decomposition of the spatial variable into radially directed blocks
of various sizes. To begin with, we denote the spatial frequency cutoff
by:
\begin{equation}
	p_{\lambda}(\xi) = \varphi_{2\lambda}(|\xi|)
	- \varphi_{\frac{1}{2}\lambda}(|\xi|) \ , \label{space_cutoff}
\end{equation}
with $P_\lambda$ the corresponding operator. For a given parameter 
$\delta\leqslant \lambda$,
we now decompose $P_\lambda$ radially as follows. First decompose the the
unit sphere $\ \ S^{n-1} \subset \mathbb{R}^n\ \ $ into angular sectors of size
$\ \ \frac{\delta}{\lambda}\times \ldots 
\times \frac{\delta}{\lambda}\ \ $ with
bounded overlap (independent of $\delta$).  
These angular sectors are then projected out to frequency $\lambda$ via rays 
through the origin. The result is a decomposition of $supp \{p_\lambda\}$
into radially directed blocks of size 
$\ \ \lambda\times\delta\times
\ldots \times \delta\ \ $ with bounded overlap. We enumerate these
blocks and label the corresponding partition of unity by 
$b^\omega_{\lambda,\delta}$. It is clear that things may be
arranged so that upon rotation
onto the $\xi_1$--axis, each $b^\omega_{\lambda,\delta}$ satisfies the 
bound:
\begin{align}\label{block_smoothness}
	|\partial_1^N b^\omega_{\lambda,\delta}| &\leqslant C_N
	\lambda^{-N} \ ,   
	& |\partial_i^N b^\omega_{\lambda,\delta}| &\leqslant C_N
	\delta^{-N} \ .
\end{align}
In particular, each $B^\omega_{\lambda,\delta}$ is given
by convolution with an $L^1$ kernel.
We shall also denote: 
\begin{align}
	S^\omega_{\lambda,d} &= B^\omega_{\lambda,(\lambda d)^\frac{1}{2}} 
	S_{\lambda,d} \ , 
	& S^\omega_{\lambda,\bullet \leqslant d} &= 
	B^\omega_{\lambda,(\lambda d)^\frac{1}{2}} 
	S_{\lambda,\bullet \leqslant d} \ . \notag
\end{align} 
Note that the operators  $S^\omega_{\lambda,d}$ and 
$S^\omega_{\lambda,\bullet \leqslant d}$ are 
only supported in the region where $\ \ |\tau| \lesssim |\xi|\ \ $.\\

We now use these multipliers to define the following dyadic norms, which will
be the building blocks for the function spaces we will use here. 

\begin{align}
	\lp{u}{X^{\frac{1}{2}}_{\lambda,p}}^p &=
	\sum_{ d\in 2^\mathbb{Z}} d^{\frac{p}{2}}
	\lp{S_{\lambda,d} u}{L^2}^p \ , &\hbox{(``classical'' $H^{s,\delta}$)}
	\label{Xsd_norm} \\ \notag \\
	\lp{u}{Y_\lambda} &=  
	\lambda^{-1}\lp{\Box S_\lambda u}{L^1(L^2)} \ ,
	&\hbox{(Duhamel)} \label{duhamel_norm} \\ \notag \\
	\lp{u}{Z_\lambda} &= \lambda^\frac{2-n}{2} \sum_d 
	\left(\sum_\omega
	\lp{S^\omega_{\lambda,d}u}{L^1(L^\infty)}^2\right)^\frac{1}{2} \ .
	&\hbox{(outer block)} \label{nullc_norm}
\end{align}
Notice that the (semi) norms $X^{\frac{1}{2}}_{\lambda,p}$ and $Y_\lambda$ 
are only well defined modulo measures supported on the light cone
in Fourier space. Because of this, it will be convenient for us
to include an extra $L^\infty(L^2)$ norm in the definition of  our 
function spaces. This represents the inclusion in the above norms of
solutions to the wave equation with $L^2$ initial data. Adding everything 
together, we are led to define the following fixed frequency (semi) norms:
\begin{equation}
	\lp{u}{F_\lambda} = \left( X^\frac{1}{2}_{\lambda,1}
	+ Y_\lambda \right) \cap S_\lambda
	\left( L^\infty(L^2) \right) \ . \label{fixedF_norm} 
\end{equation}
Unfortunately, the above norm is still not strong enough for us to be
able to iterate equations of the form \eqref{generic_system} which
contain derivatives. This is due to a very specific $Low\times High$
frequency interaction in quadratic non--linearities. 
Fortunately, this problem has been effectively handled by Tataru in
\cite{Tataru}, based on ideas from \cite{KMQij} and \cite{KTat}. 
What is necessary is to add some extra 
$L^1(L^\infty)$ norms on ``outer block'' regions of Fourier space.
This is the essence of the norm \eqref{nullc_norm} above, which is a slight
variant of that which appeared in \cite{Tataru}. This
leads to our second main dyadic norm:
\begin{equation}
	\lp{u}{G_\lambda} = \left( X^\frac{1}{2}_{\lambda,1}
	+ Y_\lambda \right) \cap S_\lambda
	\left( L^\infty(L^2) \right)\cap Z_\lambda \ . \label{fixedG_norm} 
\end{equation}
Finally, the spaces we will iterate in are produced by adding the appropriate
number of derivatives combined with the necessary Besov structures:
\begin{align}
	\lp{u}{F^s}^2 \ &= \ \sum_\lambda \lambda^{2s} \lp{u}{F_\lambda}^2
	\ , \label{F_norm} \\
	\lp{u}{G^s} \ &= \ \sum_\lambda \lambda^{s} \lp{u}{G_\lambda}
	\ . \label{G_norm}
\end{align}\\

Due to the need for precise microlocal decompositions, of crucial
importance to us will be the boundedness of certain 
multipliers on the components \eqref{Xsd_norm}--\eqref{duhamel_norm}
of our function spaces as well as mixed Lebesgue spaces. 
We state these as follows:\\

\begin{lem}[Multiplier boundedness]\label{mult_lemma}
\ \ \ 
\begin{enumerate}
	\item The following multipliers are given by $L^1$ kernels: \ \ 
	$\lambda^{-1}\nabla S_\lambda$, $S^\omega_{\lambda,d}$,
	$S^\omega_{\lambda,\bullet \leqslant d}$,
	and $(\lambda d)\varXi^{-1}S^\omega_{\lambda,d}$ . In particular,
	all of these are bounded on every mixed Lebesgue space $L^q(L^r)$. 
	\item The following multipliers are bounded on the spaces
	$L^q(L^2)$, for $1 \leqslant q \leqslant \infty$: \ \ 
	$S_{\lambda,d}$ and
	$S_{\lambda,\bullet \leqslant d}$. 
\end{enumerate}
\end{lem}\ret

\begin{proof}[Proof of Lemma \ref{mult_lemma} (1)]
First, notice that after a rescaling, the symbol for the multiplier 
$\lambda^{-1}\nabla S_\lambda$ is a $C^\infty$ bump function with $O(1)$
support. Thus, its kernel is in $L^1$ with norm independent of $\lambda$.\\

For the remainder of the operators listed in (1) above, it suffices to
work with $(\lambda d)\varXi^{-1}S^\omega_{\lambda,d}$. The boundedness of
the others follows from a similar argument. We let $\chi^\pm$ denote the
symbol of this operator cut off in the upper resp. lower half plane. After
a rotation in the spatial domain, we may assume that the spatial projection
of $\chi^\pm$ is directed along the positive $\xi_1$ axis. Now
look at $\chi^+(s,\eta)$ with coordinates:
\begin{align}
	s \ &= \ \frac{1}{\sqrt{2}}(\tau - \xi_1) \ , \notag \\
	\eta_1 \ &= \ \frac{1}{\sqrt{2}}(\tau + \xi_1) \ , \notag \\
	\eta' \ &= \ \xi' \ . \notag
\end{align} 
It is apparent that $\chi^+(s,\eta)$ has support in a box of dimension 
$\sim $
$\lambda\times\sqrt{\lambda d}
\times\ldots\times\sqrt{\lambda d}\times d$ with
sides parallel to the coordinate axis and
longest side in the $\eta_1$ direction and shortest side in the $s$ direction.
Furthermore, a direction calculation shows that one has the bounds:
\begin{align}
	|\partial^N_{\eta_1} \chi^+| &\leqslant C_N\lambda^{-N} \ ,
	&|\partial^N_{\eta'} \chi^+| &\leqslant C_N(\lambda d)^{-N/2} \ ,
	&|\partial^N_{s} \chi^+| &\leqslant C_N d^{-N} \ . \notag
\end{align}
Therefore, we have that $\chi^+$ yields an $L^1$ kernel. A similar 
argument works for the cutoff function $\chi^-$, using the rotation:
\begin{align}
	s \ &= \ \frac{1}{\sqrt{2}}(\tau + \xi_1) \ , \notag \\
	\eta_1 \ &= \ \frac{1}{\sqrt{2}}(-\tau + \xi_1) \ , \notag \\
	\eta' \ &= \ \xi' \ . \notag
\end{align} 
\end{proof}\ret

\begin{proof}[Proof of Lemma \ref{mult_lemma} (2)]
We will argue here for $S_{\lambda,d}$. The estimates for the others follow
similarly. If we denote by $K^\pm(t,x)$ the convolution kernel associated with
$S^\pm_{\lambda,d}$, then a simple calculation shows that: 
\begin{equation}
	e^{\mp2\pi i t |\xi|}\widehat{K^\pm}(t,\xi) \ = \
	\int e^{2\pi i t\tau}\psi(\tau,\xi)\, d\tau \ , \notag
\end{equation}
where $supp\{\psi\}$ is contained in a box of dimension $\ \sim\ $
$\ \ \lambda\times\ldots\times\lambda\times d\ \ $ 
with sides along the coordinate
axis and short side in the $\tau$ direction. Furthermore, one has the
estimate:
\begin{equation}
	|\partial^N_\tau \psi | \ \leqslant \ C_N \, d^{-N} \ . \notag
\end{equation}
This shows that we have the bound:
\begin{equation}
	\lp{ \widehat{K^\pm} }{L^1_\tau(L^\infty_\xi)} \ \lesssim \ 1 
	\ , \notag
\end{equation}
independent of $\lambda$ and $d$. Thus, we get the desired bounds 
for the convolution kernels.
\end{proof}\ret

As an immediate application of the above lemma, we show that the extra
$Z_\lambda$ intersection in the $G_\lambda$ norm above only effects the 
$X^\frac{1}{2}_{\lambda,1}$ portion of things.\\

\begin{lem}[Outer block estimate on $Y_\lambda$]
For $5 < n$, one has the following uniform inclusion:
\begin{equation}\label{Y_outerblock_est}
	Y_\lambda \ \subseteq \ Z_\lambda \ . 
\end{equation}
\end{lem}\ret

\begin{proof}[proof of \eqref{Y_outerblock_est}]
It is enough to show that: 
\begin{equation}
	\left( \sum_\omega \lp{\varXi^{-1} S^\omega_{\lambda,d}u}
	{L^1(L^\infty)}^2
	\right)^\frac{1}{2} \ \lesssim \ \lambda^\frac{n-4}{2} 
	\left(\frac{d}{\lambda}\right)^\frac{n-5}{4}
	\lp{S_\lambda u}{L^1(L^2)} \ . \notag
\end{equation}
First, using a local Sobolev embedding, we see that:
\begin{equation}
	\lp{B^\omega_{\lambda,(\lambda d)^\frac{1}{2}}
	\varXi^{-1}S^\omega_{\lambda,d}
	u}{L^1(L^\infty)} \ \lesssim \ \lambda^{\frac{n+1}{4}}d^\frac{n-1}{4}
	\lp{\varXi^{-1}S^\omega_{\lambda,d}u}{L^1(L^2)} \ . \notag
\end{equation}
Therefore, using the boundedness Lemma \ref{mult_lemma},
it suffices to note that by Minkowski's inequality we can bound:
\begin{align}
	\left( \sum_\omega (\int \lp{S^\omega_{\lambda,d} u}{L_x^2})^2 
	\right)^\frac{1}{2} \ &\lesssim \ \int \left(
	\sum_\omega \lp{S^\omega_{\lambda,d} u}{L_x^2}^2 
	\right)^\frac{1}{2} \ , \notag \\
	&\lesssim \ \lp{S_{\lambda,d} u}{L^1(L^2)} \ . \notag
\end{align}
\end{proof}\ret

The last line of the above proof showed that it is possible to bound
a square sum over an angular decomposition of a given function in
$L^1(L^2)$. It is also clear that this same procedure works for the 
$X^\frac{1}{2}_{\lambda,1}$ spaces because one can use Minkowski's
inequality for the $\ell^1$ sum with respect to the cone variable $d$.
This fact will be of great importance in what follows and
we record it here as:\\

\begin{lem}[Angular reconstruction of norms]
Given a test function $u$ and parameter $\delta\leqslant \lambda$,
one can bound:
\begin{equation}\label{angular_reconstruction}
	\left( \sum_\omega \lp{B^\omega_{\lambda,\delta}u}
	{X^\frac{1}{2}_{\lambda,1},Y_\lambda}^2 
	\right)^\frac{1}{2} \ \lesssim \ \lp{u}
	{X^\frac{1}{2}_{\lambda,1},Y_\lambda} \ .
\end{equation}
\end{lem}\ret

\ret

\section{Structure of the $F_\lambda$ spaces}\label{structure_section}

The purpose of this section is to clarify some remarks of the previous
section and write down two integral formulas for functions in the
$F_\lambda$ space. This material is all more or less
standard in the literature and we include it here primarily because the
notation will be useful for our scattering result. Our first order of
business is to write down a decomposition for functions in the 
$F_\lambda$ space:\\

\begin{lem}[$F_\lambda$ decomposition]\label{F_decomp}
For any $u_\lambda\in F_\lambda$, one can write:
\begin{equation}\label{u_decomp}
	u_\lambda \ = 
	\ u_{\mathring{X}_\lambda} + 
	u_{X^{1/2}_{\lambda,1}} + u_{Y_\lambda} \ ,
\end{equation}
where $u_{\mathring X_\lambda}$ is a solution to 
the homogeneous wave equation, 
$u_{X^{1/2}_{\lambda,1}}$ is the Fourier transform of an $L^1$ function, 
and $u_{Y_\lambda}$ satisfies:
\begin{equation}
	u_{Y_\lambda}(0) = \partial_t u_{Y_\lambda}\, (0) = 0 \ . \notag
\end{equation}
Furthermore, one has the norm bounds:
\begin{equation}
	\frac{1}{C} \lp{u_\lambda}{F_\lambda} \ \leqslant \
	\left(\lp{u_{\mathring X_\lambda}}{L^\infty(L^2)} + 
	\lp{u_{X^{1/2}_{\lambda,1}}}{X^\frac{1}{2}_{\lambda,1}} 
	+ \lp{u_{Y_\lambda}}{Y_\lambda}\right)
	\ \leqslant \ C \lp{u_\lambda}{F_\lambda} \ .
\end{equation}
\end{lem}\ret

We now show that the two inhomogeneous terms on the right hand side
of \eqref{u_decomp} can be written as integrals over solutions to
the wave equation with $L^2$ data. This fact will be of crucial importance
to us in the sequel. The first formula is simply a restatement of 
\eqref{Duhamel_integral}:\\

\begin{lem}[Duhamel's principle]
Using the same notation as above, for any $u_{Y_\lambda}$, one can write:
\begin{equation}\label{duhamel_int}
	u_{Y_\lambda}(t) \ = \ 
	- \int_0^t |D_x|^{-1} \sin\big( (t-s)|D_x|\big)\, \Box
	u_{Y_\lambda}(s)\ ds \ . 
\end{equation}
\end{lem}\ret

\noindent
Likewise, one can write the $u_{X^{1/2}_{\lambda,1}}$ portion of the sum
\eqref{u_decomp} as an integral over modulated solutions to the 
wave equation be foliating Fourier space by forward and backward
facing light--cones:

\begin{lem}[$X^\frac{1}{2}_{\lambda,1}$ Trace lemma]
For any $u_{X^{1/2}_{\lambda,1}}$, let 
$u^\pm_{X^{1/2}_{\lambda,1}}$ denote its restriction to the 
frequency half space $0 < \pm\tau$. Then one can write:
\begin{equation}\label{trace_formula}
	u^\pm_{X^{1/2}_{\lambda,1}}(t) \ = \
	\int e^{2\pi i t (s \pm |D_x|)}\, u^\pm_{\lambda,s}\ ds \ ,
\end{equation} 
where $u^\pm_{\lambda,s}$ is the spatial Fourier transform of 
$\widetilde{u^\pm}$ restricted to the $s^{th}$ translate of the 
forward or backward light--cone light cone in Fourier space, i.e.:
\begin{equation}
	\widehat{u^\pm}_{\lambda,s}(\xi) \ = \ \int
	\delta (\tau - s \mp |\xi|)
	\, \widetilde{u^\pm}(\tau,\xi)\ d\tau \ . \notag
\end{equation}
In particular. one has the formula:
\begin{equation}\label{X_trace_bound}
	\int \lp{u^\pm_{\lambda,s}}{L^2}\, ds  \ \lesssim \ 
	\lp{ u^\pm_{X^{1/2}_{\lambda,1}} }{X^\frac{1}{2}_{\lambda,1}} \ .
\end{equation}
\end{lem}\ret

\ret

\section{Strichartz estimates}

Our inductive estimates will be based on a method of bilinear 
decompositions and local Strichartz estimates as in the work \cite{Tataru}.
We first state the standard Strichartz from which the local estimates
follow.\\

\begin{lem}[Homogeneous Strichartz estimates (see \cite{KT_str})]
Let $5 < n$, $\sigma = \frac{n-1}{2}$, and suppose $u$ is a given function
of the spatial variable only. Then if
$\frac{1}{q}  + \frac{\sigma}{r} \leqslant \frac{\sigma}{2}$ and
$\frac{1}{q} + \frac{n}{r} = \frac{n}{2} - \gamma$, the following
estimate holds:
\begin{equation}\label{str_est}
	\lp{e^{\pm2\pi i t |D_x|}\,
	P_{\bullet \leqslant \lambda} u}{L_t^q(L_x^r)} \ \lesssim \
	\lambda^\gamma \, \lp{P_{\bullet \leqslant \lambda} u }{L^2} \ . 
\end{equation}
\end{lem}\ret

\noindent
Combining the $L^2(L^\frac{2(n-1)}{n-3})$ endpoint of the above estimate
with a local Sobolev in the spatial domain, we arrive at the following local
version of  \eqref{str_est}:\\

\begin{lem}[Local Strichartz estimate]
Let $5 < n$, 
then the following estimate holds:
\begin{equation}\label{local_str_est}
	\lp{e^{\pm2\pi i t |D_x|}\,
	B^\omega_{\lambda,(\lambda d)^\frac{1}{2}} 
	u }{L_t^2(L_x^\infty)}
	\ \lesssim \
	\lambda^\frac{n+1}{4} d^\frac{n-3}{4}
	\, \lp{B^\omega_{\lambda,(\lambda d)^\frac{1}{2}} u }{L^2} \ .
\end{equation}
\end{lem}\ret\ret

\noindent
Using the integral formulas \eqref{duhamel_int} and 
\eqref{trace_formula}, we can transfer the above
estimates to the $F_\lambda$ spaces:\\ 

\begin{lem}[$F_\lambda$ Strichartz estimates]\label{Flam_str_lem}
Let $5 < n$ and set $\sigma = \frac{n-1}{2}$. Then if
$\frac{1}{q}  + \frac{\sigma}{r} \leqslant \frac{\sigma}{2}$ and
$\frac{1}{q} + \frac{n}{r} = \frac{n}{2} - \gamma$, the following
estimates hold:
\begin{align}
	\lp{S_\lambda u}{L^q(L^r)} \ &\lesssim \
	\lambda^\gamma \, \lp{u}{F_\lambda} \ , \label{F_str_est} \\
	\left(\sum_\omega
	\lp{S^\alpha_{\lambda,\bullet \leqslant d}u}{L^2(L^\infty)}^2 
	\right)^\frac{1}{2} \ &\lesssim \
	\lambda^\frac{n+1}{4} d^\frac{n-3}{4}
	\, \lp{u}{F_\lambda} \ . \label{F_local_str_est} 
\end{align}
\end{lem}\ret

\begin{proof}[Proof of Lemma \ref{Flam_str_lem}]
It suffices to prove the estimate \eqref{F_str_est}, as the estimate 
\eqref{F_local_str_est} follows from this and a local Sobolev embedding
combined with the resumming formula \eqref{angular_reconstruction}.
Using the decomposition \eqref{F_decomp}
and the angular reconstruction formula \eqref{angular_reconstruction}, 
it is enough to prove \eqref{F_str_est} for functions
$u_{X^{1/2}_{\lambda,1}}$ and $u_{Y_\lambda}$. Using the integral formula
\eqref{trace_formula}, we see immediately that:
\begin{align}
	\lp{ u_{X^{1/2}_{\lambda,1}} }{L^q(L^r)} \ \leqslant \
	&\sum_\pm \ \int \lp{ e^{\pm 2 \pi i t |D_x|}\, u^\pm_{\lambda,s}}
	{L^q(L^r)}\ ds \ , \notag \\
	\lesssim \ &\lambda^\gamma\, \sum_\pm \ \int 
	\lp{ u^\pm_{\lambda,s} }{L^2}\ ds \ , \notag \\
	\lesssim \ &\lambda^\gamma\,  
	\lp{ u_{X^{1/2}_{\lambda,1}} }{ X^\frac{1}{2}_{\lambda,1} } \ . \notag
\end{align}
For the  $u_{Y_\lambda}$ portion of things, we can chop the function up into
a fixed number of space--time angular sectors using $L^1$
convolution kernels. Doing this and using $R_\alpha$ to denote an operator
from the set $\{I,\partial_i |D_x|^{-1}\}$, we estimate:
\begin{align}
	\lp{u_{Y_\lambda}}{L^q(L^r)} \ \leqslant \ &\lambda^{-1}
	\sum_\alpha \ \lp{\partial_\alpha u_{Y_\lambda} }{L^q(L^r)} 
	\ , \notag \\
        \lesssim \ &\lambda^{-1} \sum_{\pm\, ,\, \alpha} \ \int \ 
	\lp{ e^{\pm 2 \pi i t|D_x|}\left(
	e^{\mp 2\pi i s|D_x|} R_\alpha \,
	\Box u_{Y_\lambda}\, (s,x)\right) }{L_t^q(L_x^r)}\ ds \ , \notag \\
	\lesssim \ &\lambda^\gamma\lambda^{-1} \sum_\alpha\ \int \ 
	\lp{   e^{\mp 2\pi i s|D_x|} R_\alpha\,
	\Box u_{Y_\lambda}\, (s,x) }{L^2_x} \ ds \ , \notag \\
	\leqslant \ &\lambda^\gamma\ \lp{u_{Y_\lambda}}{Y_\lambda} \ . \notag
\end{align}
\end{proof}\ret

\noindent A consequence of \eqref{F_str_est} is that we have the embedding:
\begin{equation}\label{fixed_freq_energy_est}
	X^\frac{1}{2}_{\lambda,1} \ \subseteq \
	L^\infty(L^2) \ . \notag
\end{equation}
Using a simple approximation argument along with uniform convergence, we
arrive at the following energy estimate for the $F^s$ and $G^s$ spaces:\\

\begin{lem}[Energy estimates]\label{energy_est_lem}
For space--time functions $u$, one has the following estimates:
\begin{align}
	\lp{u}{C(\dot{H}^s)\cap C^{(1)}(\dot{H}^{s-1})} \ \lesssim \
	\lp{u}{F^s} \ , \label{Fs_energy_est} \\
	\lp{u}{C(\dot{B}^s)\cap C^{(1)}(\dot{B}^{s-1})} \ \lesssim \
	\lp{u}{G^s} \ . \label{Gs_energy_est}
\end{align}
\end{lem}\ret

\noindent Also, by duality and the estimate \eqref{fixed_freq_energy_est}, 
we have that:
\begin{equation}
	\lambda \varXi^{-1} L^1(L^2) 
	\ \subseteq \ \lambda \varXi^{-1} 
	X^{-\frac{1}{2}}_{\lambda,\infty} \ \subseteq \
	X^\frac{1}{2}_{\lambda,\infty} \ .
\end{equation}
This proves shows: \\

\begin{lem}[$L^2$ estimate for $Y_\lambda$]
The following inclusion holds uniformly:
\begin{equation}\label{Y_L2_est}
	d^\frac{1}{2} S_{\lambda,d}(Y_\lambda) \ \subseteq \
	L^2(L^2) \ , 
\end{equation}
in particular, by dyadic summing one has:
\begin{equation}
	d^\frac{1}{2} S_{\lambda,d \leqslant \bullet}
	(F_\lambda)
	\ \subseteq \ L^2(L^2) \ . \notag
\end{equation}
\end{lem}\ret

\ret

\section{Scattering}

It turns out that our scattering result, Theorem \ref{scattering_result},
is implicitly contained in the function spaces $F^s$ and $G^s$. That is,
there is scattering in these spaces independently of any specific equation 
being considered. Therefore, to prove Theorem \ref{scattering_result}, it will
only be necessary to show that our solution to \eqref{generic_system} 
belongs to these spaces.\\
 
Using a simple approximation argument, it suffices to
deal with things at fixed frequency. Because the estimates
in Theorem \ref{scattering_result} deal with more than one derivative, 
we will show that:\\

\begin{lem}[$F_\lambda$ scattering]\label{fixed_freq_scattering}
For any function $u_\lambda \in F_\lambda$, there exists a set of initial
data \ $(f_\lambda^\pm,g_\lambda^\pm)\ 
\in \ P_\lambda(L^2)\times\lambda P_\lambda (L^2)$ 
\ such that the following asymptotic holds:
\begin{align}
	\lim_{t\to \infty}\lp{u_\lambda(t)- W(f_\lambda^+,g_\lambda^+)(t)}
	{\dot{H}^1\cap\partial_t(L^2)} \ &= \ 0 \ , \label{plus_scat_dyadic}\\
	\lim_{t\to -\infty}\lp{u_\lambda(t)- W(f_\lambda^-,g_\lambda^-)(t)}
	{\dot{H}^1\cap\partial_t(L^2)} \ &= \ 0 \ . \label{minus_scat_dyadic}
\end{align}
\end{lem}\ret

\begin{proof}[Proof of Lemma \ref{fixed_freq_scattering}]
Using the notation of Section \ref{structure_section}, we may write:
\begin{equation}
	u_\lambda \ = \ u_{\mathring X_\lambda} + u^+_{X^{1/2}_{\lambda,1}} +
	u^-_{X^{1/2}_{\lambda,1}} + u_{Y_\lambda} \ , \notag
\end{equation}
We now define the scattering data implicitly by the relations:
\begin{align}
	W(f_\lambda^+,g_\lambda^+)(t) 
	\ &= \ u_{\mathring X_\lambda} - \int_{0}^\infty |D_x|^{-1}
	\sin\left(|D_x|(t-s)\right)\Box u_{Y_\lambda}(s)\, ds \ , \notag \\
	W(f_\lambda^-,g_\lambda^-)(t) 
	\ &= \ u_{\mathring X_\lambda} - \int_{-\infty}^0 |D_x|^{-1}
	\sin\left(|D_x|(t-s)\right)\Box u_{Y_\lambda}(s)\, ds \ . \notag
\end{align}
Using the fact that $\Box u_{Y_\lambda}$ has finite $L^1(L^2)$ norm, 
it suffices to show that one has the limits:
\begin{equation}
	\lim_{t\to \pm\infty}
	\lp{u^+_{X^{1/2}_{\lambda,1}}(t) + u^-_{X^{1/2}_{\lambda,1}}(t)}
	{\dot{H}^1\cap\partial_t(L^2)} \ = \ 0 \ . \notag
\end{equation}
Squaring this, we see that we must show the limits:
\begin{align}
	&\lim_{t\to\pm\infty}\int |D_x| u^+_{X^{1/2}_{\lambda,1}}(t)
	\ \overline{|D_x| u^\pm_{X^{1/2}_{\lambda,1}}}(t) \ = \ 0 
	\ , \label{first_Xlim} \\
	&\lim_{t\to\pm\infty}\int \partial_t u^+_{X^{1/2}_{\lambda,1}}(t)
	\ \overline{\partial_t u^\pm_{X^{1/2}_{\lambda,1}}}(t) \ = \ 0 
	\ . \label{second_Xlim}
\end{align}
We'll only deal here with the limit \eqref{first_Xlim}, as the limit 
\eqref{second_Xlim} follows from a virtually identical argument. 
Using the trace formula \eqref{trace_formula} along with
the Plancherel theorem, we compute:
\begin{multline}
	(L.H.S.)\eqref{first_Xlim} \\ = \ \lim_{t\to\pm\infty}
	\int e^{2\pi i t(|\xi| \mp |\xi|)}\, |\xi|^2 \
	\int e^{2\pi i t s_1} \  
	\widehat{u^+_{\lambda,s_1 + s_2}}(\xi) \
	\overline{ \widehat{u^\pm_{\lambda,s_2}}(\xi) } 
	\ \ ds_1\, ds_2 \ d\xi . \notag
\end{multline}
By \eqref{X_trace_bound} we have the bounds:
\begin{equation}
	\Lp{ |\xi|^2 \ \int   
	\big| \widehat{u^+_{\lambda,s_1 + s_2}}(\xi)\big|\ \cdot \ 
	\big| \widehat{u^\pm_{\lambda,s_2}}(\xi) \big|
	\ \ ds_1\, ds_2}{L_\xi^1} \ \lesssim \ \lambda^2\
	\lp{u^+_{X^{1/2}_{\lambda,1}}}{ X^\frac{1}{2}_{\lambda,1} } \, 
	\lp{u^\pm_{X^{1/2}_{\lambda,1}}}{X^\frac{1}{2}_{\lambda,1}}
	\ . \notag 
\end{equation}
This shows that the function:
\begin{equation}
        H_t(\xi) \ = \ |\xi|^2 \ \int   e^{2\pi i t s_1} 
	\widehat{u^+_{\lambda,s_1 + s_2}}(\xi) \
	\overline{ \widehat{u^\pm_{\lambda,s_2}}(\xi) } 
	\ \ ds_1\, ds_2 \ , \notag
\end{equation}
is bounded pointwise by an $L^1$ function uniformly in $t$. Therefore, by the
dominated convergence theorem, it suffices to show that we in fact have that
$\lim_{t\to\pm\infty} H_t = 0$. To see this, notice that by the above bounds
in conjunction with Fubini's theorem, we have that the integral:
\begin{equation}
        |\xi|^2 \ \int   
	\widehat{u^+_{\lambda,s_1 + s_2}}(\xi) \
	\overline{ \widehat{u^\pm_{\lambda,s_2}}(\xi) } 
	\ \ ds_2 \ ,  \notag
\end{equation} 
is in $L^1_{s_1}$. The result now follows from the Riemann Lebesgue Lemma.
Explicitly, one has that for almost every fixed $\xi$, the following limit
holds:
\begin{equation}
	\lim_{t\to\pm\infty} H_t(\xi) \ = \ \lim_{t\to\pm\infty}\,
	|\xi|^2 \ \int   e^{2\pi i t s_1} 
	\widehat{u^+_{\lambda,s_1 + s_2}}(\xi) \
	\overline{ \widehat{u^\pm_{\lambda,s_2}}(\xi) } 
	\ \ ds_1\, ds_2 \ \ = \ \ 0 \ . \notag
\end{equation} 
\end{proof}\ret

\ret

\section{Inductive Estimates I}

Our solution to \eqref{generic_system} 
will be produced through the usual procedure
of Picard iteration. Because the initial data and our function
spaces are both invariant with respect to the scaling \eqref{scale_trans},
any iteration procedure must effectively be global in time. Therefore,
we shall have no need of an auxiliary time cutoff system as in the
works \cite{KMopt}--\cite{Tataru}. Instead, we write \eqref{generic_system} 
directly as an integral equation:
\begin{equation}\label{generic_system_int}
	\phi \ = \ W(f,g) + \Box^{-1} \mathcal{N}(\phi,D\phi) \ . 
\end{equation}
By the contraction mapping principle and the 
quadratic nature of the nonlinearity, to produce a solution to
\eqref{generic_system_int} which satisfies the regularity assumptions
of our main theorem, it suffices to prove the following two sets
of estimates:\\

\begin{thm}[Solution of the division problem]\label{main_thm}
Let $5 < n$, then the $F$ and $G$ spaces solve the division problem for 
quadratic wave equations in the sense that for any of the model systems
we have written above: YM, WM, or MD, one has the following estimates:
\begin{align}
	\lp{\Box^{-1}\mathcal{N}(u,Dv)}{G^{s_c}} \ &\lesssim \
	\lp{u}{G^{s_c}}\lp{v}{G^{s_c}} \ , \label{division_estimate1} \\
	\lp{\Box^{-1}\mathcal{N}(u,Dv)}{F^s} \ &\lesssim \
	\lp{u}{G^{s_c}}\lp{v}{F^s} + 
	\lp{u}{F^s}\lp{v}{G^{s_c}}\ . \label{division_estimate2}
\end{align}
\end{thm}\ret\ret

The remainder of the paper is devoted to the proof of Theorem
\ref{main_thm}. In what follows, we will work exclusively with the 
equation:
\begin{equation}\label{NLW_inegral}
	\phi \ = \ W(f,g) + \Box^{-1} (\phi\, \nabla \phi) \ . 
\end{equation}
In this case, we set $s_c = \frac{n-2}{2}$. The proof of Theorem
\ref{main_thm} for the other model equations can be achieved through
a straightforward adaptation of the estimates we give here. In fact, after
the various derivatives and values of $s_c$ are taken into account, the
proof in these cases follows verbatim from estimates \eqref{HH_est} and 
\eqref{HL_est} below. \\

Our first step is to take a Littlewood-Paley
decomposition of $\Box^{-1} (u \, \nabla v)$ with respect to
space--time frequencies:
\begin{equation}\label{LP_decomp1}
	\Box^{-1} (u \, \nabla v) \ = \ \sum_{\mu_i}
	\Box^{-1}(S_{\mu_1} u \, \nabla S_{\mu_2} v) \ .
\end{equation}
We now follow the standard procedure of splitting the sum
\eqref{LP_decomp1} into three pieces depending on the cases
$\mu_1 \ll \mu_2$, $\mu_2 \ll \mu_1$, and
$\mu_2 \sim \mu_1$.
Therefore, due to the $\ell^1$ Besov structure in the $F$ spaces, in order 
to prove both \eqref{division_estimate1} and 
\eqref{division_estimate2}, it suffices to show the two estimates:\\

\begin{align}
	\lp{\Box^{-1}(S_{\mu_1} u \, \nabla S_{\mu_2} v)}{G_\lambda}
	\ &\lesssim \ 
	\lambda^{-1} \mu_1^\frac{n}{2} 
	\lp{u}{F_{\mu_1}}\lp{v}{F_{\mu_2}} \ \ , \
	\mu_1 \sim \mu_2  \ ,
	\label{HH_est} \\ \notag \\ 
	\lp{\Box^{-1}(S_{\mu} u \, \nabla S_{\lambda} v)}{G_\lambda}
	\ &\lesssim \ \mu^\frac{n-2}{2}
	\lp{u}{G_{\mu}}\lp{v}{F_{\lambda}} \ \ , \
	\mu \ll \lambda \ . \label{HL_est}
\end{align}\\

\noindent Notice that after some weight trading,
the estimates \eqref{division_estimate1} and 
\eqref{division_estimate2}  follow from \eqref{HL_est} in the case where 
$\mu_2 \ll \mu_1$.\\

\begin{proof}[proof of \eqref{HH_est}]
It is enough if we show the following two estimates:

\begin{align}
	\lp{S_\lambda(S_{\mu_1} u \, \nabla S_{\mu_2} v)}{L^1(L^2)}
	\ &\lesssim \ \mu_1^\frac{n}{2} 
	\lp{u}{F_{\mu_1}}\lp{v}{F_{\mu_2}} \ \ , \
	\mu_1 \sim \mu_2  \ , 
	\label{HH_est_reduced1} \\
	\lp{S_\lambda \Box^{-1}(S_{\mu_1} u \, \nabla S_{\mu_2} v)}
	{L^\infty(L^2)}
	\ &\lesssim \ \lambda^{-1} \mu_1^\frac{n}{2} 
	\lp{u}{F_{\mu_1}}\lp{v}{F_{\mu_2}} \ \ , \
	\mu_1 \sim \mu_2  \ .
	\label{HH_est_reduced2}
\end{align}\\

\noindent
In fact, it suffices to prove \eqref{HH_est_reduced1}. To see this, notice
that one has the formula:
\begin{equation}
	\left[S_\lambda,\Box^{-1}\right] G \ = \ 
	W( E*S_\lambda G) - S_\lambda W(E*G)\ . \notag
\end{equation}
Thus, after multiplying by $S_\lambda$, we see that:
\begin{align}
	S_\lambda \left[S_\lambda,\Box^{-1}\right] G \ &= \
	P_\lambda\big( W( E*S_\lambda G) - S_\lambda W(E*G)\big) 
	\ , \notag \\
	&= \ W( E*S_\lambda P_\lambda G) - S_\lambda W(E*P_\lambda G) 
	\ , \notag \\
	&= \ S_\lambda \left[S_\lambda,\Box^{-1}\right] P_\lambda G 
	\ . \notag
\end{align}
Therefore, by the (approximate)
idempotence of $S_\lambda$ one has: 
\begin{align}\label{commutator_equation}
	S_\lambda \Box^{-1} G \ &= \
	S_\lambda \Box^{-1} S_\lambda G + S_\lambda
	\left[S_\lambda,\Box^{-1}\right] G \ , \notag \\
	&= \ S_\lambda \Box^{-1} S_\lambda G + 
	S_\lambda \left[S_\lambda,\Box^{-1}\right] P_\lambda G
	\ . \notag
\end{align}
Thus, by the boundedness of $S_\lambda$ on the spaces
$L^\infty(L^2)$ the energy estimate, one can bound:
\begin{equation}
	\lp{S_\lambda \Box^{-1} G}{L^\infty(L^2)}
	\ \lesssim \ \lambda^{-1}\left(\lp{S_\lambda G}{L^1(L^2)} 
	+ \lp{P_\lambda G}{L^1(L^2)} + \lp{S_\lambda P_\lambda G}{L^1(L^2)}\right) 
	\ . \notag
\end{equation}\\

\noindent
We now use the fact that the multipliers $S_\lambda$ and $P_\lambda$ are both bounded on
the space $L^1(L^2)$ to reduce things to the estimate:
\begin{align}
	\lp{S_{\mu_1} u \, \nabla S_{\mu_2} v}{L^1(L^2)} \ \lesssim \ \
	&\lp{S_{\mu_1} u}{L^2(L^4)}\lp{\nabla S_{\mu_2} v}{L^2(L^4)} \ ,
	\notag \\
	& \ \mu_1^\frac{n-2}{4}
	\mu_2^\frac{n+2}{4} \, \lp{u}{F_{\mu_1}}\lp{v}{F_{\mu_2}}
	\ . \notag 
\end{align}
Taking into account the the bound $\mu_1 \sim \mu_2$, the 
claim now follows.
\end{proof}\ret

Next, we'll deal with the estimate \eqref{HL_est}. For the remainder of the 
paper we shall fix both $\lambda$ and $\mu$ and assume they such that
$\mu \ll \lambda$ for a fixed constant. We now decompose the product
$S_\lambda(S_\mu u \, \nabla S_\lambda v)$ into a sum of three pieces:
\begin{equation}
	S_\lambda(S_\mu u \, \nabla S_\lambda v) \ = \ A + B + C \ , \notag
\end{equation}
where
\begin{align}
	A \ &= \ S_\lambda(S_\mu u \, \nabla 
	S_{\lambda,c\mu \leqslant \bullet} v) 
	\ , \notag \\ \notag \\
	B \ &= \ S_{\lambda,c\mu \leqslant \bullet}(
	S_\mu u \, \nabla S_{\lambda,\bullet < c\mu } v) 
	\ , \notag \\ \notag \\	
	C \ &= \ S_{\lambda,\bullet < c\mu}(
	S_\mu u \, \nabla S_{\lambda,\bullet < c\mu} v) \ . \notag
\end{align}
Here $c$ is a suitably small constant which will be chosen later. It will be
needed to make explicit a dependency between some of the constants which
arise in a specific frequency localization in the sequel.
We now work to recover the estimate \eqref{HL_est} for 
each of the three above terms separately.

\begin{proof}[proof of \eqref{HL_est} for the term $A$]
Following the remarks at the beginning of the proof of \eqref{HH_est},
it suffices to compute:
\begin{align}
	\lp{S_\lambda(S_\mu u \, \nabla 
	S_{\lambda,c\mu \leqslant \bullet} v)}{L^1(L^2)} \ &\lesssim \
	\lambda \ \lp{S_\mu u}{L^2(L^\infty)} 
	\lp{S_{\lambda,c\mu \leqslant \bullet} v}{L^2(L^2)} \ , \notag \\
	&\lesssim \ \lambda \mu^\frac{n-1}{2}\lp{u}{F_\mu}\
	(c\mu)^{-\frac{1}{2}}\lp{v}{F_\lambda} \ , \notag \\
	&\lesssim \
	c^{-1} \lambda \mu^\frac{n-2}{2} \lp{u}{F_\mu}\lp{v}{F_\lambda} 
	\ . \notag 
\end{align}
For a fixed $c$, we obtain the desired result.
\end{proof}\ret

We now move on to showing the inclusion \eqref{HL_est} for the $B$ term above.
In this range, we are forced to work outside the context of $L^1(L^2)$
estimates. This is the reason we have included the $L^2(L^2)$ based
$X^\frac{1}{2}_{\lambda,1}$ spaces. This also means that we will need to 
recover $Z_\lambda$ norms by hand (because they are only covered by the
$Y_\lambda$ spaces). However, because this last task will require a somewhat
finer analysis than what we will do in this section, we contend ourselves here
with showing:\\

\begin{proof}[proof of the $X^\frac{1}{2}_{\lambda,1}
\cap S_\lambda(L^\infty(L^2))$ estimates for the term $B$]
Our first task will be deal with the energy estimate
which we write as:
\begin{equation}
	\lp{S_\lambda \Box^{-1} S_{\lambda,c\mu \leqslant \bullet}(
	S_\mu u \, \nabla S_{\lambda,\bullet < c\mu } v) }{L^\infty(L^2)}
	\ \lesssim \ \mu^\frac{n-2}{2} \lp{u}{F_\mu} \lp{v}{F_\lambda}
	\ . \notag
\end{equation}
For $G$ supported away from the light--cone in Fourier space, we have
the identity:
\begin{equation}
	S_\lambda \Box^{-1}S_\lambda G \ 
	= \ \varXi^{-1}S_\lambda G - W(\varXi^{-1}P_\lambda S_\lambda 
	G) \ . \notag
\end{equation}
Therefore, by using the energy estimate for the $X^\frac{1}{2}_{\lambda,1}$ space, 
this allows us to estimate: 
\begin{align}
	\lp{S_\lambda \Box^{-1} S_\lambda G}{L^\infty(L^2)} \ &\lesssim \
	\lp{\varXi^{-1}S_\lambda G}{L^\infty(L^2)} + \lp{W(\varXi^{-1}P_\lambda S_\lambda 
	G)}{L^\infty(L^2)} \ , \notag \\
	&\lesssim \ \lp{\varXi^{-1}S_\lambda G}{L^\infty(L^2)} \ , \notag \\
	&\lesssim \ \lp{\varXi^{-1}S_\lambda G}{X^\frac{1}{2}_{\lambda,1}}\ . \notag
\end{align}
Therefore, we are left with estimating the
term $B$ in the $X^\frac{1}{2}_{\lambda,1}$
space. For a fixed distance $d$ from the cone, we compute that:
\begin{align}
	\lp{\varXi^{-1}S_{\lambda,d} S_{\lambda,c\mu \leqslant \bullet}(
	S_\mu u \, \nabla S_{\lambda,\bullet < c\mu } v) }{L^2(L^2)} 
	\ &\lesssim \
	d^{-1} \ \lp{S_\mu u}{L^2(L^\infty)}\lp{S_\lambda v}{L^\infty(L^2)}
	\ , \notag \\
	&\lesssim \ d^{-1} \mu^\frac{n-2}{2} 
	\lp{v}{F_\mu} \lp{u}{F_\lambda} \ . \notag
\end{align}
Summing $d^\frac{1}{2}$ times this last 
expression over all $c\mu \leqslant d$ yields:
\begin{multline}
	\sum_{c\mu \leqslant d} d^\frac{1}{2}
	\lp{\varXi^{-1}S_{\lambda,d} S_{\lambda,c\mu \leqslant \bullet}(
	S_\mu u \, \nabla S_{\lambda,\bullet < c\mu } v) }{L^2(L^2)}
	\\ \lesssim \ \sum_{c\mu \leqslant d}
	\left( \frac{\mu}{d}\right)^\frac{1}{2} \mu^\frac{n-2}{2}
	\lp{v}{F_\mu} \lp{u}{F_\lambda} \ . \notag
\end{multline}
For a fixed $c$ we obtain the desired result.\\
\end{proof}\ret

\ret

\section{Interlude: Some bilinear decompositions}

To proceed further, it will be necessary for us to take a closer look at the 
expression:
\begin{align}\label{B_leftovers}
	&S^\omega_{\lambda,d}(
	S_\mu u \, \nabla S_{\lambda,\bullet < c\mu } v) \ , 
	&c\mu \leqslant d \ ,
\end{align}
as well as the sum:
\begin{equation}\label{Cdecomp}
	C \ = \ S_{\lambda,\bullet < c\mu}(
	S_\mu u \, \nabla S_{\lambda,\bullet < c\mu} v) \ = \
	C_{I} + C_{II} + C_{III} \ , \notag
\end{equation}
where
\begin{align}
	C_{I} &= \sum_{d < c\mu}
	S_{\lambda,d}(S_{\mu,\bullet \leqslant d}u \ \nabla
	S_{\lambda, \bullet \leqslant d} v) \ , \notag \\ \notag \\
	C_{II} &= \sum_{d < c\mu}
	S_{\lambda, \bullet < d}(
	S_{\mu,\bullet \leqslant d}u \ \nabla
	S_{\lambda,d} v) \ , \notag \\ \notag \\
	C_{III} &= \sum_{d\leqslant \mu}
	S_{\lambda,\bullet < min\{c\mu,d\}}(
	S_{\mu,d} u \ \nabla S_{\lambda,\bullet < min\{c\mu,d\}}
	v ) \ . \notag
\end{align}\\

We'll begin with a decomposition of $C_I$ and $C_{II}$. The $C_{III}$ term is
basically the same but requires a slightly more delicate analysis. All of
the decompositions we compute here will be for a fixed $d$. The full 
decomposition will then be given by summing over the relevant values of $d$.
Because our decompositions will
be with respect to Fourier supports, it suffices to look at the convolution
product of the corresponding cutoff functions in Fourier space. In what 
follows, we'll only deal with the $C_I$ term. It will become apparent that the
same idea works for $C_{II}$. Therefore,
without loss of generality, we shall decompose the product:
\begin{equation}\label{mult_prod1}
	s^+_{\lambda,d}( s^\pm_{\mu,\bullet \leqslant d} * 
	s^+_{\lambda,\bullet \leqslant d}) \ . 
\end{equation}
To do this, we use the standard device of restricting the angle of interaction 
in the above product. It will be crucial for us to be able to make
these restrictions based only on the spatial Fourier variables, because we
will need to reconstruct our decompositions through square--summing.
For $\ \ (\tau',\xi')\in 
supp\{s^\pm_{\mu,\bullet \leqslant d}\}\ \ $ and $\ \ (\tau,\xi)\in
supp\{s^+_{\lambda,\bullet \leqslant d}\}\ \ $ we compute that:

\begin{align}
	O(d) &= \big| | \tau' + \tau | - |\xi' + \xi| \big| \ , 
	\notag \\ 
	&= \Big| \big| \pm|\xi'| + |\xi| + O(d) \big|
	- |\xi' + \xi| \Big| \ , \notag \\
	&= \Big| O(d) + \big| \pm |\xi'| + |\xi| \big| - |\xi' + \xi| \Big| 
	\ . \notag 
\end{align}\\

\noindent Using now the fact that $\ \ d < c\mu\ \ $ and 
$\ \ \mu < c\lambda\ \ $ to conclude that $\ \ |\xi'| \sim \mu\ \ $ and
$\ \ |\xi| \sim \lambda\ \ $, we see that one has the angular restriction:
\begin{equation}
	\mu \Theta^2_{\pm\xi',\xi} \ \lesssim \
	\Big| \pm |\xi'| + |\xi|  - |\xi' + \xi|\Big|  \ = \ O(d) \ . \notag
\end{equation}
In particular we have that $\ \ \Theta_{\pm\xi',\xi} 
\lesssim \sqrt{\frac{d}{\mu}}\ \ $. This allows us to decompose the product 
\eqref{mult_prod1} into a sum over angular regions with
$\ \ O(\sqrt{\frac{d}{\mu}})\ \ $  spread. The result is:\\

\begin{lem}[Wide angle decomposition]\label{main_decomp_lem1}
In the ranges stated for the $C_I$ term above, one can write:
\begin{equation}\label{main_decom1}
	s^+_{\lambda,d}( s^\pm_{\mu,\bullet \leqslant d} * 
	s^+_{\lambda,\bullet \leqslant d}) \ = \
	\sum_{\substack{ \omega_1,\omega_2,\omega_3 \ : \\
	|\omega_1\mp\omega_2|\sim (d/\mu)^\frac{1}{2} \\
	|\omega_1 - \omega_3| \sim (d/\mu)^\frac{1}{2} } }
	b^{\omega_1}_{\lambda,\lambda (\frac{d}{\mu})^\frac{1}{2}}
	s^+_{\lambda,d}\left(
	{s^{\omega_2}_{\mu,\bullet \leqslant d}}^\pm \ * \ 
	b^{\omega_3}_{\lambda,\lambda (\frac{d}{\mu})^\frac{1}{2}}
	s^+_{\lambda,
	\bullet \leqslant d}\right) \ . 
\end{equation}\\
for the convolution of the associated cutoff functions in Fourier space.
\end{lem}\ret

\noindent 
We note here that the key feature in the decomposition \eqref{main_decom1} is 
that the sum is (essentially) diagonal in all three angles which
appear there ($\omega_1,\omega_2,\omega_3$).
It is useful here to keep in mind the following diagram:\\

\begin{figure}[h]
        \scalebox{.80}{\hspace{.4in}\includegraphics{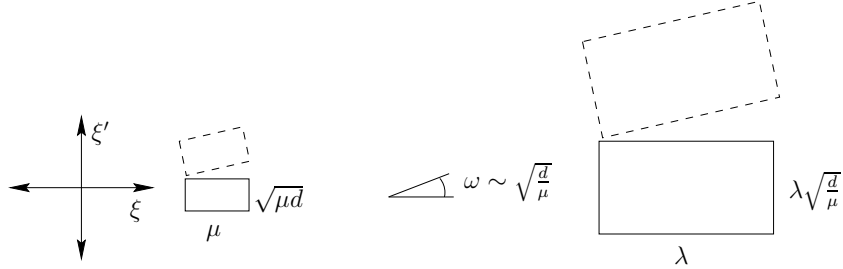}}
        \caption{Spatial supports in the wide angle decomposition.}
	\label{wide_angle_decomp}
\end{figure}\ret

We now focus our attention on decomposing the convolution:
\begin{equation}\label{mult_prod2}
	s^+_{\lambda,\bullet \leqslant min\{c\mu,d\}}(
	s^\pm_{\mu,d} * s^+_{\lambda,\bullet \leqslant min\{c\mu,d\}}) \ .
\end{equation}
If it is the case that $\ \ d \ll \mu\ \ $, then the 
same calculation which was used to
produce \eqref{main_decom1} works and we end up with the same type of 
sum. However, if we are in the case where $\ \ d\sim\mu\ \ $, 
we need to compute things a bit more carefully in order to ensure that 
we may still decompose the multiplier $s^\pm_{\mu,d}$ using only restrictions
in the spatial variable. To do this, we will now assume that things are set up so that
$\ \ c\mu \ll d\ \ $. It is clear that 
all the previous decompositions can be made
so that we can reduce things to this consideration.
If we now take $\ \ (\tau',\xi')\in 
supp\{s^\pm_{\mu,d}\}\ \ $ and $\ \ (\tau,\xi)\in
supp\{s^+_{\lambda,\bullet \leqslant min\{c\mu,d\}}\}\ \ $, we can use the
facts that $\ \ \tau' = O^\mp(d) \pm |\xi'| \ \ $, 
$\ \ \tau = O(c\mu) + |\xi|\ \ $, and
$\ \ |\xi| \gg \mu\ \ $ to compute that:\\

\begin{align}
	O(c\mu) &= \Big|  \big| \tau' + \tau 
	\big|  -  |\xi' + \xi| \Big| \ , \notag \\
	&= \Big|  O^\mp(d) 
	\pm|\xi'| + O(c\mu) + |\xi|   -  |\xi' + \xi| \Big| \ , 
	\label{Opm_force_line}
\end{align}\\

\noindent 
where the term $O^\mp(d)$ in the above expression is such that
$\ \ |O^\mp(d)| \sim d\ \ $. In fact, one can see that the equality
\eqref{Opm_force_line} forces $\ \ \pm O^\mp(d) < 0\ \ $ on account of the fact
that $\ \ \pm (\pm|\xi'| + |\xi| -  |\xi' + \xi|) > 0\ \ $ and the assumption
$\ \ |O(c\mu) + O^\mp(d)|\sim d\ \ $. In particular,
this means that we can multiply $s^\pm_{\mu,d}$
in the product \eqref{mult_prod2} by the cutoff $s_{|\tau| < |\xi|}$ 
without effecting things. This in turn shows that we may decompose the product
\eqref{mult_prod2} based solely on restriction of the spatial Fourier
variables, just as we did to get the sum in Lemma
\ref{main_decomp_lem1}. \\

We now return to the $C_I$ term. For the sequel, we will need to know what 
the contribution of the factor $S_{\lambda, \bullet \leqslant d}v$ to the
following localized product is:
\begin{equation}
	S^{\omega_1}_{\lambda,d}(S_{\mu,\bullet \leqslant d}u \ \nabla
	S_{\lambda, \bullet \leqslant d} v) \ . \notag
\end{equation}
Using Lemma \ref{main_decomp_lem1}, we see that we may write:
\begin{equation}\label{first_try}
	s^{\omega_1}_{\lambda,d}(s_{\mu,\bullet \leqslant d}*
	s_{\lambda, \bullet \leqslant d}) \ = \
	s^{\omega_1}_{\lambda,d}\left(
	s^{\omega_2}_{\mu,\bullet \leqslant d} \ * \ 
	b^{\omega_3}_{\lambda,\lambda (\frac{d}{\mu})^\frac{1}{2}}
	s_{\lambda,\bullet \leqslant d}\right) \ ,
\end{equation}
where $\ \ |\omega_1 - \pm\omega_2| \sim |\omega_3 - \pm\omega_2| \sim 
\sqrt{\frac{d}{\mu}}\ \ $.  However, this can be refined significantly. To
see this, assume that the spatial support of $s^{\omega_1}_{\lambda,d}$
lies along the positive $\xi_1$ axis. We'll label this block by 
$b^{\omega_1}_{\lambda,(\lambda d)^\frac{1}{2}}$. Because we are in the range where
$\ \ \sqrt{\mu d} \ll \sqrt{\lambda d}\ \ $, we see that since for any 
$\ \ \xi\in supp\{b^{\omega_3}_{\lambda,\lambda (\frac{d}{\mu})^\frac{1}{2}}
\}\ \ $ and $\ \ \xi'\in supp_{\xi'}\{ s^{\omega_2}_{\mu,\bullet \leqslant d} 
\}\ \ $ the sum $\ \ \xi + \xi'\ \ $
must belong to $\ \  supp\{b^{\omega_1}_{\lambda,(\lambda d)^\frac{1}{2}}\}\ \ $, we in fact
have that $\xi$ itself must belong to a block of size $\ \ \lambda\times
\sqrt{\lambda d}\times\ldots\times\sqrt{\lambda d} \ \ $.
This allows us to write: \\

\begin{lem}[Small angle decomposition]\label{main_decomp_lem2}
In the ranges stated for the $C_I$ term above, we can write:
\begin{equation}\label{main_decomp2}
	s^{\omega_1}_{\lambda,d}(s_{\mu,\bullet \leqslant d} *
	s_{\lambda, \bullet \leqslant d} ) \ = \
	s^{\omega_1}_{\lambda,d}(
	s^{\omega_2}_{\mu,\bullet \leqslant d} \ * \ 
	s^{\omega_3}_{\lambda,\bullet \leqslant d}) \ .
\end{equation}
where $\ \ |\omega_1 - \omega_3|\sim \sqrt{\frac{d}{\lambda}}\ \ $, and
$\ \ |\omega_1 - \pm\omega_2|\sim \sqrt{\frac{d}{\mu}}\ \ $.
\end{lem}\ret

\noindent
It is important to note here that if one were to sum the expression
\eqref{main_decomp2} over $\omega_1$, the resulting sum would be 
(essentially) diagonal in $\omega_3$, but there would be \emph{many}
$\omega_1$ which would contribute to a single $\omega_2$. This means that
the resulting would \emph{not} be diagonal in $\omega_2$ as was the case
for the sum \eqref{main_decom1}. It is helpful to visualize things through the 
following figure:\\

\begin{figure}[h]
        \scalebox{.80}{\hspace{.2in}\includegraphics{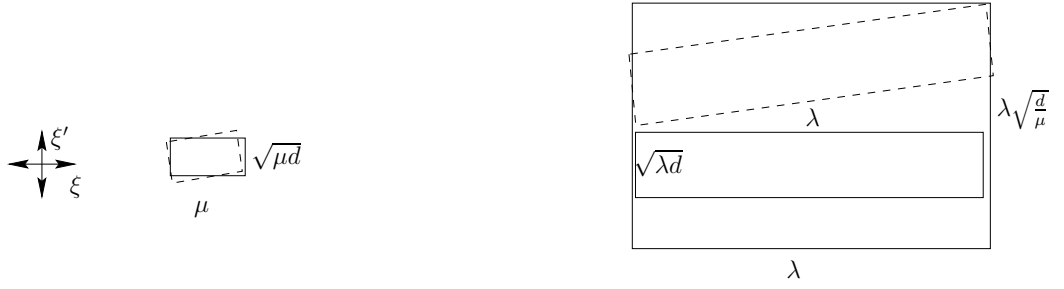}}
        \caption{Spatial supports in the small angular decomposition.}
	\label{angle_decomp}
\end{figure}\ret

Our final task here is to mention an analog of Lemma \ref{main_decomp_lem2}
for the term \eqref{B_leftovers}. Here we can frequency localize the factor
$S_{\lambda,\bullet < c\mu}$ in the product using the fact that one
has $\ \ \mu \ll c^{-\frac{1}{2}}\sqrt{\lambda d}\ \ $. The result is:\\

\begin{lem}[Small angle decomposition for the term $B$]\label{main_decomp_lem3}
In the ranges stated for the $B$ term above, we can write:
\begin{equation}\label{main_decomp3}
	s^{\omega_1}_{\lambda,d}(s_{\mu} *
	s_{\lambda, \bullet \leqslant c\mu} ) \ = \
	s^{\omega_1}_{\lambda,d}(
	s_{\mu} \ * \ b^{\omega_3}_{\lambda, (\lambda d)^\frac{1}{2}}
	s_{\lambda,\bullet \leqslant c\mu}) \ .
\end{equation}
where $\ \ |\omega_1 - \omega_3|\sim \sqrt{\frac{d}{\lambda}}\ \ $. 
\end{lem}\ret

\noindent 
Finally, we note here the important fact that in the decomposition
\eqref{main_decomp3} above, the range of interaction in the product forces
$\ \ d \lesssim \mu \ \ $. 
This completes our list of bilinear decompositions.

\ret

\section{Inductive Estimates II: Remainder of the
$Low \times High \Rightarrow High$ frequency interaction}

It remains for us is to bound the term $B$ from line \eqref{B_leftovers}
in the $Z_\lambda$ space,
as well as show the inclusion \eqref{HL_est} for the terms $C_I$
-- $C_{III}$ from line \eqref{Cdecomp}. 
We do this now, proceeding in reverse order. \\

\begin{proof}[proof of estimate \eqref{HL_est} for the $C_{III}$ term]
To begin with we fix $d$. Using the remarks at the beginning of the
proof of \eqref{HH_est}, we see that it is enough to show that:
\begin{multline}\label{CIIIest}
	\lp{ S_{\lambda,\bullet < min\{c\mu,d\}}(
	S_{\mu,d} u \ \nabla S_{\lambda,\bullet < min\{c\mu,d\}}
	v )}{L^1(L^2)} \\ \lesssim \ \lambda
	\left(\sum_\omega \lp{S^\omega_{\mu,d} u}
	{L^1(L^\infty)}^2\right)^\frac{1}{2}
	\ \lp{v}{F_\lambda} \ . 
\end{multline}
To accomplish this, we first use the wide angle decomposition, 
\eqref{main_decom1}, on the left hand side of \eqref{CIIIest}. 
This allows us to compute, using a Cauchy--Schwartz, that:
\begin{align}
	&\lp{ S_{\lambda,\bullet < min\{c\mu,d\}}(
	S_{\mu,d} u \ \nabla S_{\lambda,\bullet < min\{c\mu,d\}}
	v )}{L^1(L^2)}  \ , \notag \\
	\lesssim \ \ &\sum_{\substack{\omega_2 , \omega_3 \ : \\
	|\omega_3 \pm\omega_2| 
	\sim (d/\mu)^\frac{1}{2} }  } 
	\lp{S^{\omega_2}_{\mu,d} u}{L^1(L^\infty)} \cdot
	\lp{\nabla B^{\omega_3}_{\lambda,\lambda (\frac{d}{\mu})^\frac{1}{2}}
	S_{\lambda,\bullet < min\{c\mu,d\}}
	v )}{L^\infty(L^2)} \ , \notag \\
	\lesssim \ \ &\lambda
	\left(\sum_\omega \lp{S^\omega_{\mu,d} u}
	{L^1(L^\infty)}^2\right)^\frac{1}{2} \left(
	\sum_\omega 
	\lp{B^{\omega}_{\lambda,\lambda (\frac{d}{\mu})^\frac{1}{2}}
	S_{\lambda,\bullet < min\{c\mu,d\}}v )}{L^\infty(L_x^2)}^2
	\right)^\frac{1}{2} \ , \notag \\
	\lesssim \ \ &\lambda
	\left(\sum_\omega \lp{S^\omega_{\mu,d} u}
	{L^1(L^\infty)}^2\right)^\frac{1}{2}
	\lp{v}{F_\lambda} \ . \notag
\end{align}
Summing over $d$ now yields the desired estimate.
\end{proof}\ret

\begin{proof}[proof of \eqref{HL_est} for the $C_{II}$ term]
Again, fixing $d$, and using the angular decomposition lemma 
\ref{main_decomp_lem1}, we compute that:
\begin{align}
	&\lp{S_{\lambda, \bullet < d}(
	S_{\mu,\bullet \leqslant d}u \ \nabla
	S_{\lambda,d} v)}{L^1(L^2)} \ , \notag \\
	\lesssim \ 
	&\lambda \sum_{\substack{ \omega_2,\omega_3 \ : \\
	|\omega_3\pm\omega_2| \sim (d/\mu)^\frac{1}{2} }  } 
	\lp{S^{\omega_2}_{\mu,\bullet \leqslant d} u}{L^2(L^\infty)} \cdot
	\lp{B^{\omega_3}_{\lambda,\lambda (\frac{d}{\mu})^\frac{1}{2}}
	S_{\lambda,d}
	v }{L^2(L^2)} \ , \notag \\
	\lesssim \ &\lambda \left(\sum_\omega
	\lp{S^{\omega}_{\mu,\bullet \leqslant d} u}
	{L^2(L^\infty)}^2\right)^\frac{1}{2}
	\lp{S_{\lambda,d} v}{L^2(L^2)} \ , \notag \\
	\lesssim \ &\lambda\mu^\frac{n-2}{2}
	\left(\frac{d}{\mu}\right)^\frac{n-5}{4}
	\lp{u}{F_\mu}\lp{v}{F_\lambda} \ . \notag
\end{align}
This last expression can now be summed over $d$, using the condition
$d < c\mu$, to obtain the desired result.
\end{proof}\ret


\begin{proof}[proof of \eqref{HL_est} for the $C_{I}$ term]
This is the other instance where we will have to rely on the
$X^\frac{1}{2}_{\lambda,1}$ space. Following the same reasoning used
previously, we first bound:
\begin{align}
	\ \ \ &\lp{\varXi^{-1}
	S_{\lambda,d}(S_{\mu,\bullet \leqslant d}u \ \nabla
	S_{\lambda, \bullet \leqslant d} v)}{L^2(L^2)} \ , \notag \\
	\lesssim \ \ &d^{-1}    
	\sum_{\substack{ \omega_2,\omega_3 \ : \\
	|\omega_3\pm\omega_2| \sim (d/\mu)^\frac{1}{2}  } } 
	\lp{S^{\omega_2}_{\mu,\bullet \leqslant d} u}{L^2(L^\infty)} \cdot
	\lp{B^{\omega_3}_{\lambda,\lambda (\frac{d}{\mu})^\frac{1}{2}}
	S_{\lambda,\bullet \leqslant d}
	v )}{L^\infty(L^2)} \ , \notag \\
	\lesssim \ \ &d^{-1}
	\left(\sum_\omega
	\lp{S^{\omega}_{\mu,\bullet \leqslant d} u}
	{L^2(L^\infty)}^2\right)^\frac{1}{2} \cdot 
	\left( \sum_\omega 
	\lp{B^{\omega}_{\lambda,\lambda (\frac{d}{\mu})^\frac{1}{2}}
	S_{\lambda,\bullet \leqslant d}v )}{L^\infty(L_x^2)}^2 
	\right)^\frac{1}{2} \ , \notag \\
	\lesssim \ \ &d^{-\frac{1}{2}} \mu^\frac{n-2}{2}
	\left(\frac{d}{\mu}\right)^\frac{n-5}{4}
	\lp{u}{F_\mu}\lp{v}{F_\lambda} \ . \notag
\end{align}
Multiplying this last expression by $d^\frac{1}{2}$ and then
using the condition $d < c\mu$ to sum over $d$ yields the desired result
for the $X^\frac{1}{2}_{\lambda,1}$ space part of estimate \eqref{HL_est}.
It remains to prove the $Z_\lambda$ estimate. Here we use the
second angular decomposition lemma \ref{main_decomp_lem2} to compute that
for fixed $d$:
\begin{align}
	\ \ \ &\left( \sum_{\omega_1} \lp{\varXi^{-1} S^{\omega_1}_{\lambda,d}
	(S_{\mu,\bullet \leqslant d}u \ \nabla
	S_{\lambda, \bullet \leqslant d} v)}{L^1(L^\infty)}^2
	\right)^\frac{1}{2} \ , \notag \\
	\lesssim \ \ &(\lambda d)^{-1} \left( 
	\sum_{\substack{ \omega_1,\omega_2,\omega_3 \ : \\
	\omega_1 - \omega_3 \sim (d/\lambda)^\frac{1}{2} \\
	\omega_1 \pm \omega_2 \sim (d/\mu)^\frac{1}{2}} } \lp{
	S^{\omega_1}_{\lambda,d}
	(S^{\omega_2}_{\mu,\bullet \leqslant d}u \ \nabla
	S^{\omega_3}_{\lambda, \bullet \leqslant d} v)}{L^1(L^\infty)}^2
	\right)^\frac{1}{2} \ , \notag \\
	\lesssim \ \ &d^{-1} \sup_\omega\lp{
	S^{\omega}_{\mu,\bullet \leqslant d}u}{L^2(L^\infty)}\cdot
	\left( \sum_\omega \lp{
	S^{\omega}_{\lambda, \bullet \leqslant d} v}{L^2(L^\infty)}^2
	\right)^\frac{1}{2} \ , \notag \\
	\lesssim \ \ &\left(\frac{d}{\mu}\right)^\frac{n-5}{4}
	\left(\frac{d}{\lambda}\right)^\frac{n-5}{4}
	\mu^\frac{n-2}{2}\lambda^\frac{n-2}{2}
	\lp{u}{F_\mu}\lp{v}{F_\lambda} \ . \notag
\end{align}
Multiplying this last expression by $\lambda^\frac{2-n}{2}$ and
summing over $d$ using the condition
$\ \ d < \lambda,\mu\ \ $ yields the desired result.
\end{proof}\ret


\begin{proof}[proof of the $Z_\lambda$ embedding for the $B$ term]
The pattern here follows that of the last few lines of the previous proof.
Fixing $d$, we use the decomposition
Lemma \ref{main_decomp_lem3} to compute that:
\begin{align}
	\ \ \ &\left( \sum_{\omega} \lp{\varXi^{-1} S^{\omega}_{\lambda,d}
	(S_{\mu}u \ \nabla
	S_{\lambda, \bullet \leqslant c\mu} v)}{L^1(L^\infty)}^2
	\right)^\frac{1}{2} \ , \notag \\
	\lesssim \ \ &(\lambda d)^{-1} \left( 
	\sum_{\substack{ \omega_1 , \omega_3 \ : \\
	|\omega_1 - \omega_3| \sim (d/\lambda)^\frac{1}{2} }} \lp{
	S^{\omega_1}_{\lambda,d}
	(S_{\mu}u \ \nabla
	B^{\omega_3}_{\lambda,(\lambda d)^\frac{1}{2}}
	S_{\lambda, \bullet \leqslant c\mu} v)}{L^1(L^\infty)}^2
	\right)^\frac{1}{2} \ , \notag \\
	\lesssim \ \ &d^{-1} \lp{
	S_{\mu}u}{L^2(L^\infty)}\cdot
	\left( \sum_\omega \lp{B^{\omega}_{\lambda,(\lambda d)^\frac{1}{2}}
	S_{\lambda, \bullet \leqslant c\mu} v}{L^2(L^\infty)}^2
	\right)^\frac{1}{2} \ , \notag \\
	\lesssim \ \ &\left(\frac{\mu}{d}\right)^\frac{1}{2}
	\left(\frac{d}{\lambda}\right)^\frac{n-5}{4}
	\mu^\frac{n-2}{2}\lambda^\frac{n-2}{2}
	\lp{u}{F_\mu}\lp{v}{F_\lambda} \ . \notag
\end{align}
Multiplying the last line above by a factor of $\lambda^\frac{2-n}{2}$
and using the conditions $\ \ d < \lambda\ \ $ and
$\ \ c\mu < d \lesssim \mu\ \ $, 
we may sum over $d$ to yield the desired result.
\end{proof}

\ret


\end{document}